\documentclass[%
     aip,
      jmp,%
       amsmath,amssymb,
       preprint,%
    ]{revtex4-1}

\usepackage{amsmath,amsthm,amssymb}
\usepackage{booktabs}
\usepackage[english]{babel}
\usepackage{amsmath}
\usepackage{graphicx}
\usepackage[colorlinks=true, allcolors=blue]{hyperref}
\newcommand{\bq} {\mathbf{q}}

\newcommand{\bu} {\mathbf{u}}
\newcommand{\buu} {u}

\newcommand{\btau} {\mbox{\boldmath $\tau$}}
\newcommand{\bgam} {\mbox{\boldmath $\gamma$}}

\newcommand{\blam} {\mbox{\boldmath $\lambda$}}

\begin{document}
\title{Application of machine learning to viscoplastic flow modeling}
\author{E. Muravleva}
 \author{I. Oseledets}%
  \author{D. Koroteev}
%  \email{\{e.muravleva, i.oseledets, d.koroteev\}@skoltech.ru}
  \affiliation{ 
      Skolkovo Institute of Science and Technology.
  }%

   \date{\today}% It is always \today, today,

%\author{Ekaterina Muravleva, Ivan Oseledets, Dmitry Koroteev}
\begin{abstract}
We present a method to construct reduced-order models for
duct flows of Bingham media. Our method is based on proper orthogonal decomposition (POD)
to find a low-dimensional approximation to the velocity and artificial neural network
to approximate the coefficients of a given solution in the constructed POD basis.
We use well-established augmented Lagrangian method and
finite-element discretization in the ``offline'' stage.
We show that the resulting approximation has a reasonable accuracy,
but the evaluation of the approximate solution several orders of magnitude times faster.

%Numerical simulation of the flow of viscoplastic media is not a straightforward task.
\end{abstract}

\maketitle

\section{Introduction}
%Yield stress fluid flows play an important role in many applications, especially in the oil and gas industry \cite{FrigPasoMend2017}.
Viscoplastic properties of materials play an important role for various technological processes related to development of hydrocarbon reservoirs \cite{FrigPasoMend2017}. Examples of such processes are hydraulic fracturing \cite{guo2017numerical, shojaei2014continuum, osiptsov2017fluid} and flow diverters  \cite{bunger2017four,li2012viscoplastic}. 
Real-time control in such processes is often required, and it is impossible without fast solvers. 

Viscoplastic flow modeling is often computationally expensive, see recent reviews \cite{SarWachs2017,MitTsam2017}.
In this paper we
apply reduced order modeling and machine learning techniques to approximate the results of numerical simulations
of viscoplastic media. The idea is that the physical system is described by a few number of input parameters
(i.e. Bingham number), but the numerical simulation is computationally expensive.
Our goal is to compute the approximation to the solution by learning
from the results of numerical simulations for different values of parameters,
describing the system.

The computation is split into two steps.
At the first step (which is also called \emph{offline stage}) we conduct numerical simulation for different values of parameters and collect
solutions (so-called \emph{snapshots}).
This step is computationally expensive and is typically
done using high-performance computing. Based on the result of numerical simulation we
compute the approximant, which can be computed efficiently
for any new value of parameters in the \emph{online stage}.
Although this is a standard framework in the field
of reduced order modeling, application of these methods to numerical simulation of viscoplastic media faces
several challenges which we address.
The main challenge is that the equations are nonlinear, and even
it is quite straightforward to reduce the dimensionality of the problem using
\emph{proper orthogonal decomposition} (POD), it is not simple to construct the approximant that is easy to evaluate
for new parameter values.
In order to solve this problem, we introduce an additional approximant, which is based on
artificial neural networks (ANN) and is learned using standard backpropagation techniques. As an alternative to ANN
other approximation schemes for multivariate functions can be used. One of the promising approaches is the
\emph{proper generalized decomposition}, which was succesfully applied to different computational rheology problems
\cite{ammar2006new, chinesta2011overview}.
The main difference with our approach is that the parametric dependencies in the considered problems are typically smooth,
whereas ANN can approximate more general classes of functions, which are, for example, piecewise-smooth.
To summarize, main contributions of our paper are:
\begin{itemize}
\item We propose a ``black-box'' approach to construct a reduced-order model for Bingham fluid flow in a channel based
on POD to compute a low-dimensional projection, and ANN to approximate the coefficients of POD decomposition from the parameters.
\item We show the accuracy of the proposed approximation for a single yield stress for different domains.
\item We show the accuracy of the proposed approximation for piecewise-constant yield stress limit.
\item Our method allows to achieve several orders of magnitude speedup for the considered examples.
\end{itemize}

\section{Test problem}

\subsection{Governing equations}
The constitutive relations of viscoplastic Bingham medium
connect the stress tensor $\btau$ to the rate-of-strain tensor $\dot{\bgam}$ as follows:
\begin{equation*}
 \begin{split}
  & \dot{\bgam} = 0, \quad |\tau| \leq \tau_s, \\
  & \btau = \left(\frac{\tau_s}{|\dot{\gamma}|} + \mu\right) \dot{\bgam}, \quad |\tau| > \tau_s,
  \end{split}
\end{equation*}
where $\mu$ is the plastic viscosity, $\tau_s$ is the yield stress, $\bu$ is the velocity vector, $\dot{\bgam} = \frac{1}{2} \left(\nabla \buu + (\nabla \buu)^{\top}\right)$,  and the norms of the tensors $\bgam$ and $\btau$ are defined by
$$
|\dot{\bgam}| = \sqrt{\dot{\bgam} : \dot{\bgam}},\quad
|\btau| = \sqrt{\btau : \btau}.
$$
As a test problem, we consider well-known Mosolov problem  \cite{MosMyas1965, MosMyas1966, MosMyas1967}. It describes
an isothermal steady laminar flow of an incompressible Bingham fluid in an infinitely long duct 
with a cross-section $\Omega \subset \mathbb{R}^2$ under the action of the pressure gradient. It is modeled by the following equation
\begin{equation}\label{bingml:mosolov0}
    -\mu \nabla^2 u - \tau_s \nabla \cdot \left( \frac{\nabla u}{\vert \nabla u \vert} \right) = \textcolor{black}{\nabla p},
\quad u_{\partial \Omega} = 0,
\end{equation}
which follows from mass and momentum conservation laws.
We consider the no-slip (Dirichlet) boundary conditions.
Here $u$ is the axial velocity ($\bu = (0,0,\buu)$),
$\textcolor{black}{\nabla p}$ is a pressure drop, $\textcolor{black}{\nabla u = \left(\frac{\partial u}{\partial x}, \frac{\partial u}{\partial y} \right)}$.

After proper rescaling \cite{Huilgol2005} $\textcolor{black}{ \widehat{x} = \frac{x}{L}, \quad \widehat{u} = \frac{L^2 \nabla p}{\mu} 
u, \quad B = \frac{\tau_s}{L \nabla p}}$, \textcolor{black}{$L$ is a characteristic length 
of the cross-section},
we can obtain the dimensionless form
of the problem \eqref{bingml:mosolov0} with the Bingham number $B$ as a parameter:
\begin{equation}\label{bingml:mosolov}
    \textcolor{black}{
        - \nabla^2 \widehat{u} - B \nabla \cdot \left( \frac{\nabla \widehat{u}}{\vert \nabla \widehat{u} \vert} \right) = 1,
    \quad \widehat{u}_{\partial \Omega} = 0.}
\end{equation}
From now on we will use the form \eqref{bingml:mosolov} and omit caps in the notation.

This problem has been solved numerically in many papers (both for steady
\cite{SaramRoquet2001,MoyFrig2004,Huilgol2005,MuravMM2009,walton1991axial,szabo1992flow,wachs2007numerical} 
and unsteady \cite{muravleva2009unsteady,Mur_JNNFM2010} cases). Wall slip boundary conditions were considered in
\cite{roquet2008adaptive,DamPhil,DamGeor}. 

\subsection{Numerical methodology}
The solution of the Mosolov problem \textcolor{black}{\eqref{bingml:mosolov}} can be found from the minimization of the functional
\begin{equation}\label{bingml:functional}
    J(\buu) =  \int \left(\frac{1}{2}\vert \nabla \buu \vert^2 + B \vert \nabla \buu \vert - 
     \buu \right) dx \rightarrow \min_{\buu \in H^1_0 (\Omega)}
\end{equation}
One of the most widely-used approaches to solve \eqref{bingml:functional} is the augmented Lagrangian method (ALM) \cite{Glowinski}, which has the following form. First, we introduce an
additional variable $\bq = \nabla \buu$ and consider \textcolor{black}{constrained} minimization problem
$$
   \min_{\buu, \bq = \nabla u} J(u, \bq), \quad J(\buu, \bq) =
   \int \left(\frac{1}{2}\vert \bq \vert^2 + B \vert \bq \vert -  \buu \right) dx.
$$
To deal with  the constraint, the Lagrange multiplier $\blam$ is introduced, and also a penalty term is added,
thus leading to the augmented Lagrangian
\begin{equation}\label{bingml:auglag}
    \mathcal{L}_r(\buu, \bq, \blam) =  \int \left(\frac{1}{2}\vert \bq \vert^2 + B \vert \bq \vert - \buu \right) dx + 
    \int \blam \cdot (\bq - \nabla \buu) dx + \frac{r}{2} \int \vert \bq - \nabla \buu \vert^2 dx.
\end{equation}
The original minimization problem \eqref{bingml:functional} is equivalent to finding a saddle point of the functional \eqref{bingml:auglag}
$$
\max_{\blam} \min_{\buu, \bq} \mathcal{L}_r(\buu, \bq, \blam).
$$
For finding the saddle point of \eqref{bingml:auglag} we use the ALG2 method \cite{Glowinski}, which is an iterative algorithm.
At each iteration, given $\buu^{(k)}, \bq^k, \blam^k$ we compute the next approximation 
$\buu^{(k+1)}$, $\bq^{(k+1)}$, $\blam^{(k+1)}$
by the following steps:
\begin{itemize}
\item (Update $u$) Solve

\begin{equation}\label{bingml:alg2-1}
    - r \nabla^2 \buu^{(k+1)} =  1 + \nabla \cdot (\blam^k - r \bq^k).
\end{equation}
\item (Update $\bq$)
\begin{equation}\label{bingml:alg2-2}
  \bq^{(k+1)} = \begin{cases}
  0,  \mbox{if} \quad \vert \blam^k + r \nabla \buu^{(k+1)} \vert  \leq B, \\
  \frac{\blam^k + r \nabla \buu^{(k+1)}}{1+r} \left( 1 - \frac{B}{\left\vert \blam^k + 
  r \nabla \buu^{(k+1)}\right\vert}\right),  \quad  \mbox{otherwise.}
    \end{cases}
\end{equation}
  \item (Update  $\blam$)
  \begin{equation}\label{bingml:alg2-3}
     \blam^{(k+1)} = \blam^{(k)} + r \left(\bq^{(k+1)} - \nabla \buu^{(k+1)}\right).
  \end{equation}
\end{itemize}
As an initial approximation, we choose $\bq^{(0)} = 0$, $\blam^{(0)}= 0$.
To implement \eqref{bingml:alg2-1}, \eqref{bingml:alg2-2}, \eqref{bingml:alg2-3} numerically, we use weak formulations for each of these problems and finite element method (FEM). We use FENICS package \cite{LoggMardalEtAl2012a} for the implementation, where only weak formulation of the problem is necessary, and everything else (including unstructured mesh generation) is done automatically. We used standard finite element spaces: for velocity we used piecewise-linear functions, and for $\bq$ and $\blam$ -- piecewise-constant functions.  Given the solver (which consists of executing steps \eqref{bingml:alg2-1}, \eqref{bingml:alg2-2}, \eqref{bingml:alg2-3}  until convergence), we can now focus on the main problem considered in this paper: the construction of the reduced model.

\section{Construction of the reduced model}
ALG2 algorithm is time consuming, since it typically requires many iterations, and on each iteration a boundary value problem needs to be solved.
We will refer to the ALG2 method as the \emph{full model}.
Our goal is to construct a \emph{reduced model} which takes the same input parameters as the input, and computes the approximation to the solution obtained 
by the full model, but much faster. In order to construct the reduced order model we use a standard procedure. We first run the full model 
for a certain set of input parameters in the
so-called \emph{offline stage}, and then, using these results, construct the reduced order model that is able to compute approximation to the solution 
for any other parameter $B$ in the \emph{online stage}. In this paper we are focusing on the reduced order model for the computation of the velocity $\buu$.
\subsection{Dimensionality reduction}
Let us describe the proposed approximation procedure more formally. The full model can be considered as a computation of the mapping $B \rightarrow \buu(B)$,
which maps
from the space of parameters to the space of solutions $\buu(B) \in \mathbb{R}^N$, where $N$ is the number of finite elements used.
Our approximation procedure consists in two steps. First, we construct a low-dimensional subspace in $\mathbb{R}^N$ that approximates $\buu(B)$ for all $B$ of interest.
Specifically, we look for the basis vectors $u_1, \ldots, u_m$ of length $N$, where $m \ll N$  and organize them into an $N \times m$ matrix
$$
  U = \begin{bmatrix} u_1 & \ldots & u_m \end{bmatrix}.
$$
Then, for any $B$ we approximate $\buu(B)$ as
\begin{equation}\label{bingml:podapprox}
\buu(B) \approx U c(B), \quad c(B) = \begin{bmatrix} c_1(B) & \ldots & c_m(B) \end{bmatrix}.
\end{equation}
where $B \rightarrow c(B)$ is a mapping from the parameter space to $\mathbb{R}^m$. 
Representation \eqref{bingml:podapprox} can be also written as
$$
\buu(B) \approx \sum_{k=1}^m u_{k} c_{k}(B).
$$

We need to solve two problems: construct the basis and compute coefficients for given $B$ in a fast way.

The construction of such low-dimensional basis is called \emph{linear dimensionality reduction} \cite{kunisch2002galerkin, peterson1989reduced,lassila2014model},
and the representation \eqref{bingml:podapprox} is called \emph{proper orthogonal decomposition} (POD). POD is defined by the matrix $U$.
The columns of $U$ are called \emph{POD modes}. 
This matrix can be computed by using singular value decomposition of the so-called \emph{snapshot matrix}. To construct this matrix, we select
$M > m$ values of parameters $B_1, \ldots, B_M$, compute the solutions using the full model $\buu(B_1) \ldots \buu(B_M)$ and put 
them as columns into an $N \times M$ snapshot matrix $S$.
The POD modes are obtained as the first $m$  left singular vectors of this matrix. The number of modes $m$ can be estimated from the decay of singular values of $S$.

However, the dimensionality reduction is not enough, since in the online stage the coefficients $c(B)$ have to be computed. Thus, the 
second step is the approximation of $c(B)$ by 
another function $\widehat{c}(B)$ which we will be able to evaluate in a fast way. If such approximation is given, the computation of \eqref{bingml:podapprox} 
is done by first evaluating the vector $\widehat{c}(B)$ and then computing the matrix-by-vector product $U \widehat{c}(B)$.

To compute $\widehat{c}(B)$ we can reuse the snapshots which were used for the construction of the basis. 
Given $\buu(B_1), \ldots, \buu(B_M)$ and basis $U$ we compute
$$
    c(B_k) = U^{\top} \buu(B_k), \quad k = 1, \ldots, M. 
$$
Then we fix a set of mappings $B \rightarrow \widehat{c}(B)$ and find the one that minimizes the mean squared error: 
\begin{equation}\label{bingml:minfunc}
   \sum_{k=1}^M \left \Vert c (B_k) - \widehat{c}(B_k) \right \Vert^2 \rightarrow \min,
   \end{equation}
   where $\Vert \cdot \Vert$ is the Euclidean norm of the vector.
   As an approximation class of mappings $\widehat{c}(B)$, we will use artificial neural networks, 
   which are very efficient for the approximation of multivariate functions.
 
\subsection{Approximation of coefficients using artificial neural network}
A detailed introduction of artificial neural networks is out of the scope of this paper, and can be found, for example, in \cite{Goodfellow-et-al-2016},
but for the convenience we will provide basic definitions. 
Artificial neural network is a parametrized mapping from the input to the output, given as a superposition of simple functions. 
%A detailed introduction
%can be found, for example, in \cite{Goodfellow-et-al-2016}.
We use so-called deep feedforward fully-connected networks, which are defined as follows. Given an input vector $z$ of length $n_1$, 
we introduce a weight matrix $W_1$ of size $n_2 \times n_1$, a bias vector $b_1$ of length $n_2$, select a nonlinear function $f$ that maps from 
$\mathbb{R}$ to $\mathbb{R}$ and compute
$$ y_i  = f\left(b_i + \sum_{j=1}^{n_1} W_{ij} z_j\right), \quad i = 1, \ldots, n_2.$$
Function $f$ is a function of one variable and is called \emph{activation function}.
Typical choice for this function is a sigmoidal function, or so-called \emph{rectified linear unit} (ReLU). 
For the simplicity we will use notation
$$ y = f(W z+b),$$
for such transformation, meaning that $f$ is applied elementwise to the vector $W z + b$.
This is a non-linear mapping, and corresponds to a one-layer neural network. A deep neural network is defined by taking a superposition of such transforms.
We introduce matrices $W_k$ of size $n_{k+1} \times n_{k}$ and vectors $b_k$ of length $n_k, \quad k =1, \ldots, K$. Using them, we define a $K$-layer neural network transformation 
from an input vector $z$ to an output vector $y$ as
$$
    y_1 = f(W_1 z + b_1), \quad y_2 = f(W_2 y_1 + b_2),  \quad \ldots, \quad y = y_K = f(W_K y_{K-1} + b_K).
$$
The vectors $y_k, \quad k=1, \ldots, K-1$ are outputs of \emph{hidden layers} of a neural network. The size of the vector $y_k$ is called 
the width of the $k$-th layer. 
It is also very important that neural networks provide universal approximation: with sufficient width of the layers a good approximation can be obtained for any continious function \cite{hornik1989multilayer}. In practice, however, we are interested in smaller widths of the layers while maintaining the required accuracy.
We first fix the neural network structure, i.e. define the number of hidden layers and their widths, and also select the function $f$. 
The weight matrices $W_k$ and vectors $b_k$, $k =1, \ldots, K$ are called \emph{parameters} of a neural network model.

\begin{figure}[h!]
    \includegraphics[width=0.9\textwidth]{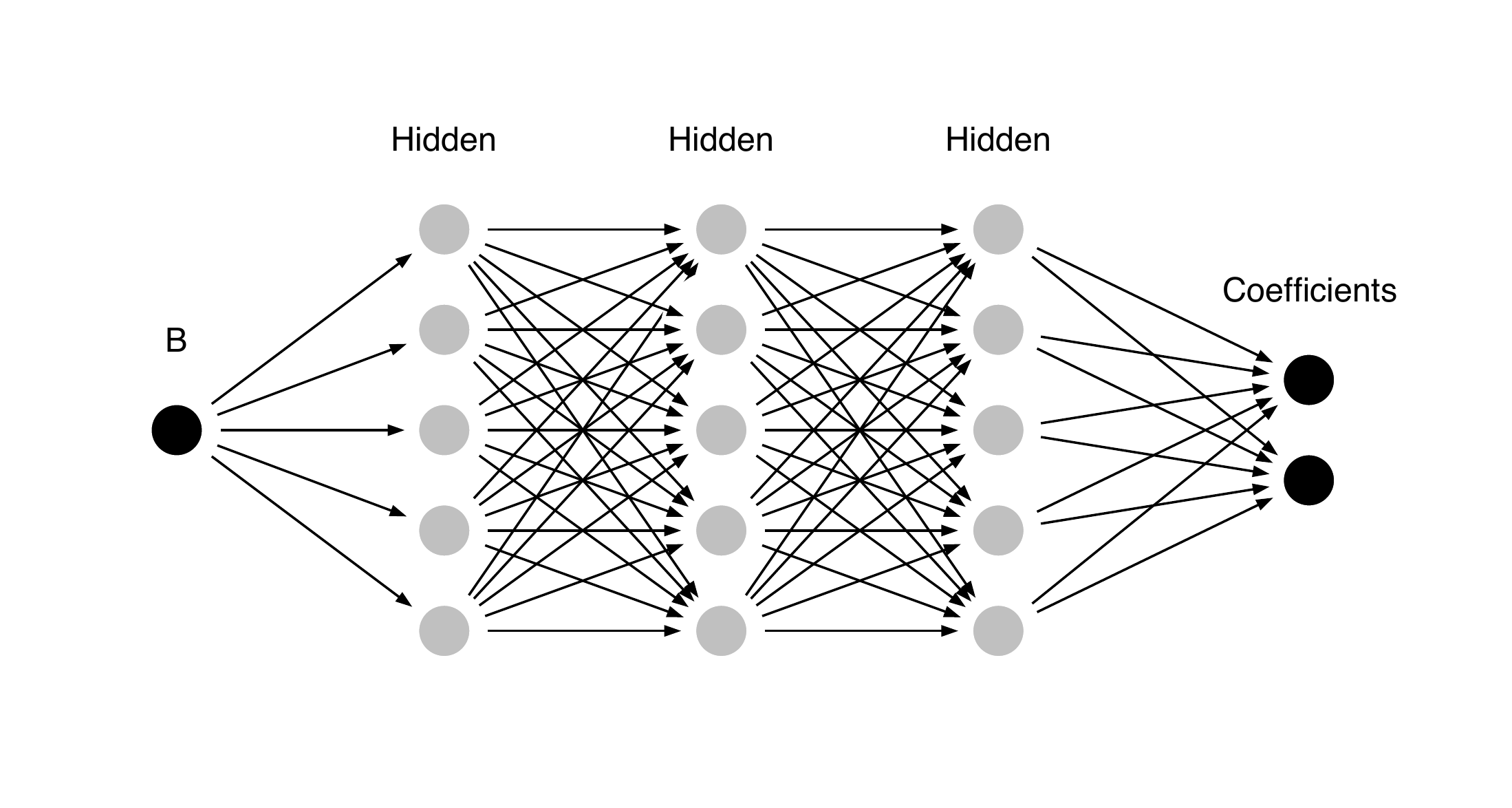}
    \caption{Three-layer neural network with one-dimensional input, hidden layer widths $5$, $5$, $5$ and output of length $2$.}\label{bingml:nnpic}
\end{figure}

In our case the neural network takes $B$ as an input, and outputs a vector of coefficients of length $m$.  
An example network  structure is shown on Figure~\ref{bingml:nnpic}. 
Given the pairs $(B_k, c(B_k)), k = 1, \ldots, M$, we minimize the functional \eqref{bingml:minfunc} with respect to the parameters of the neural network. 
This is not a simple task, since it is a non-convex optimization problem. Fortunately, there are several 
standard methods for solving this optimization problem, based on stochastic gradient methods \cite{Goodfellow-et-al-2016}, which are implemented in
modern machine learning packages, like Tensorflow \cite{tensorflow2015-whitepaper}, and can be readily used. 
These methods have very good performance in practice. Among them, ADAM optimizer \cite{kingma2014adam} is one of the most efficient and we will use it.

\section{Numerical results}
\subsection{Single yield stress limit}
\begin{figure}[h!]
\includegraphics[width=0.9\textwidth]{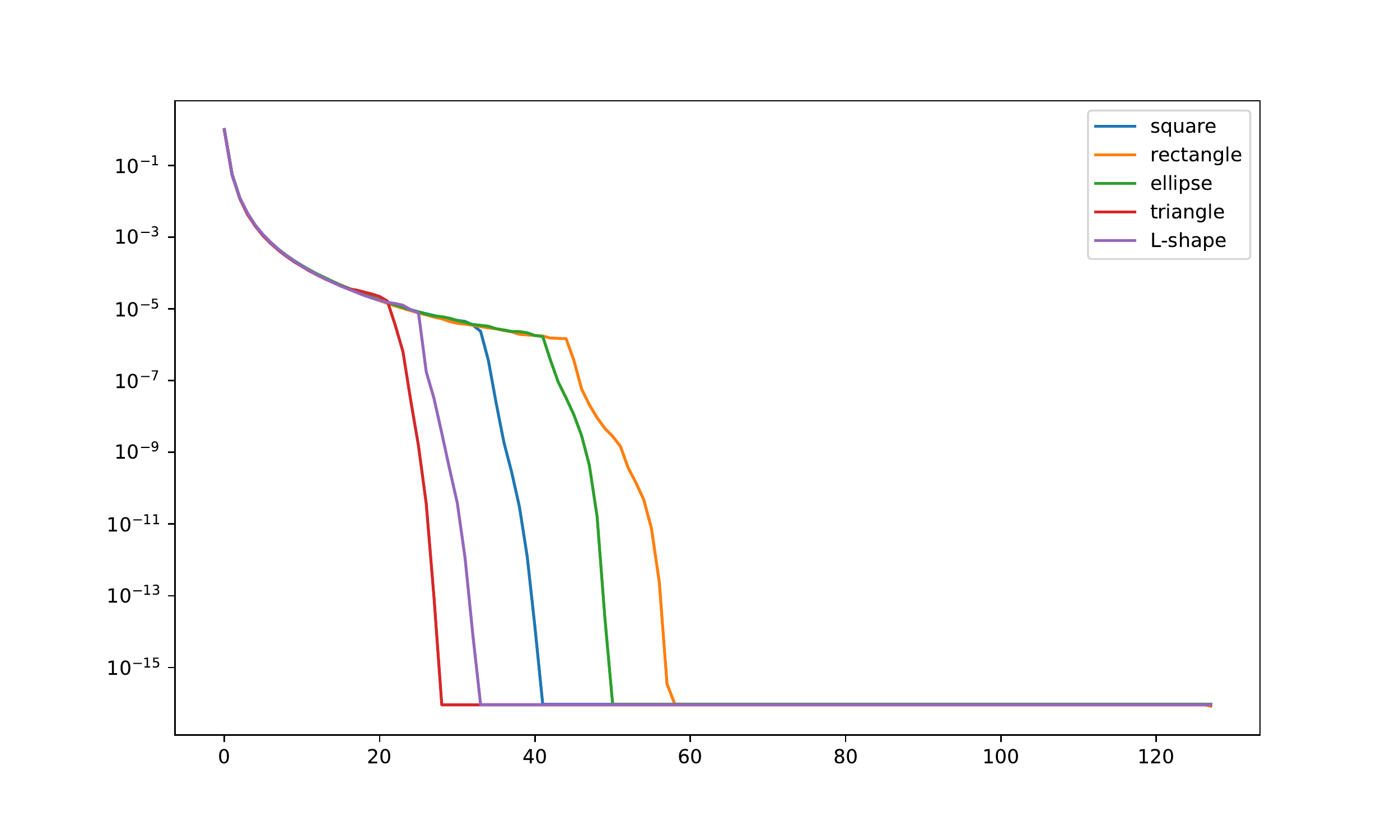}
\caption{Decay of singular values of the snapshot matrices for different domains. Horizontal axis corresponds to the number of the singular value. }\label{bingml:svdomains}.
\end{figure}

\begin{figure}[h!]
\includegraphics[width=0.9\textwidth]{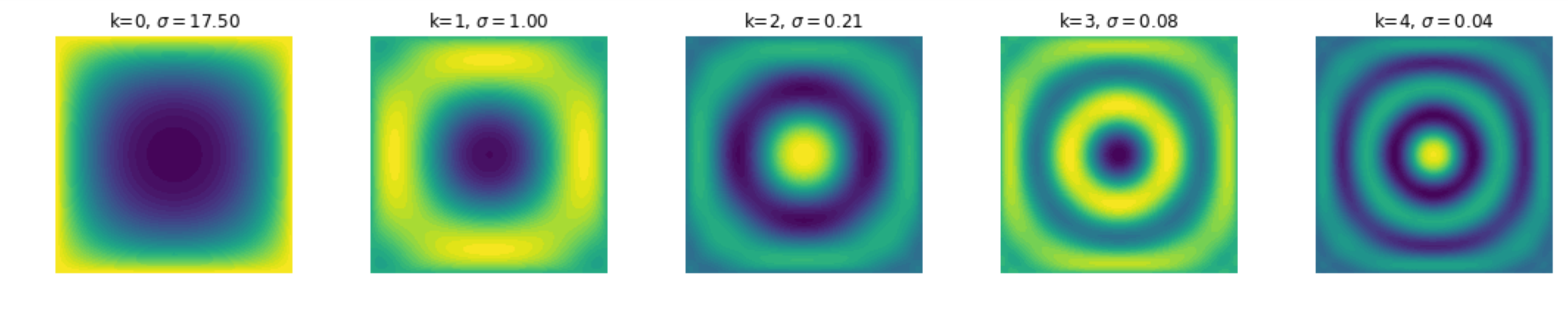} (a)
\vskip 2mm
\includegraphics[width=0.9\textwidth]{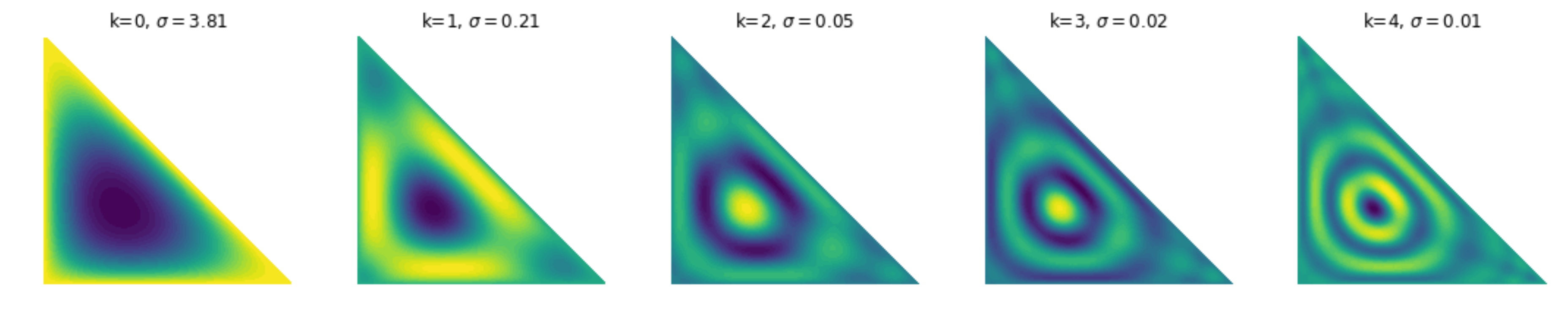} (b)
\vskip 2mm
\includegraphics[width=0.9\textwidth]{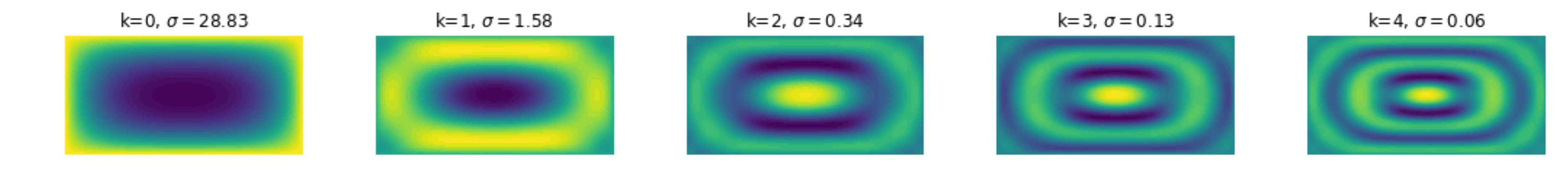}(c)
\vskip 2mm
\includegraphics[width=0.9\textwidth]{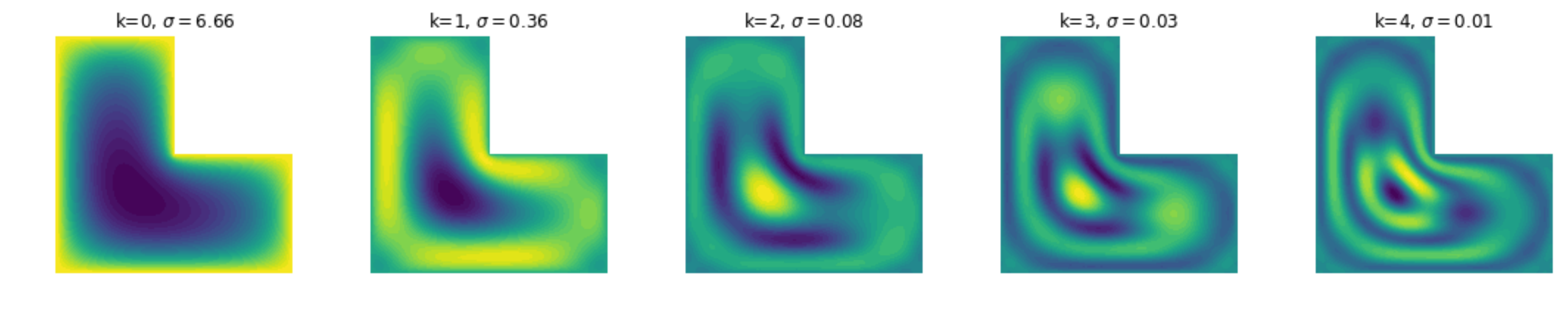} (d)
\vskip 2mm
\includegraphics[width=0.9\textwidth]{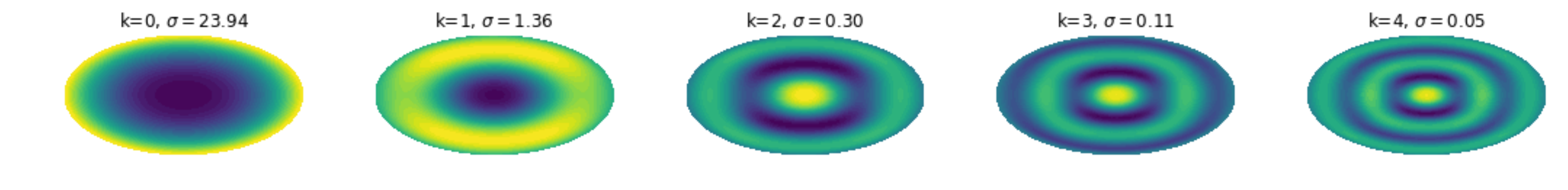} (e)
\caption{First 5 POD basis functions for different domains}\label{bingml:domains}
\end{figure}

As our first example, we consider a single-parameter case, namely, the parameter is the Bingham number $B$.
We examine several types of cross-sections: square, rectangle, triangle, L-shaped domain.
We randomly sample Bingham numbers from a uniform distribution on $[0, 1]$ to obtain snapshots. Note that for all these domains the critical yield stress is smaller than $1$, 
so for $B > B_{crit}$ the velocity is equal to $0$, and snapshots, corresponding to these parameters do not provide any additional information. Such behavior is automatically captured by the proposed algorithm. For each domain, we collect the snapshot matrices, compute first POD basis functions, and corresponding coefficients. 
Singular values of the snapshot matrices decay very fast (see Figure~\ref{bingml:svdomains}), and $m = 20$ provides very high accuracy of the approximation.
The first basis functions for different domains are shown on Figure~\ref{bingml:domains}.

As we see, the singular values decay very fast, thus it is sufficient to leave at most $m=20$
coefficients to get a good approximation. It is also interesting to look at the POD modes.
The first $5$ POD modes and corresponding singular values are shown in Figure~\ref{bingml:domains}. Note that
the first POD function has the structure ``similar'' to the typical flow pattern for each shape. This is an interesting topic for further research.

To construct a neural network approximation for the mapping from $B$ to the first $K$ coefficients of the POD decomposition, we randomly split the dataset into test and train parts (33\% test and 67\% train), and learn a corresponding mapping. As an architecture of ANN we take $3$ hidden layers with $60, 50, 40$ hidden units, respectively and one output layer. As an optimizer, ADAM optimizer \cite{kingma2014adam} was used with default parameters.
The relative error of approximation of coefficients on the test set for different domains is given in Table~\ref{bingml:errdom}.
%This error is defined as a sum over all test samples, divided by the sum of all coefficients:
%$$\varepsilon_{test} = \sqrt{\frac{\sum_{i=1}^{N_{test}} \sum_{k=1}^K \left(c^{(i)}_k - \widehat{c}^{(i)}_k \right)^2}{
%\sum_{i=1}^{N_{test}} \sum_{k=1}^K (c^{(i)}_k)^2}}.
%$$

\begin{table}[h!]
\begin{center}
\begin{tabular}{lr}
\toprule
    Domain &  Relative error \\
\midrule
    square &        0.006686 \\
  triangle &        0.005700 \\
   L-shape &        0.008667 \\
 rectangle &        0.005344 \\
   ellipse &        0.003489 \\
\bottomrule
\end{tabular}
\caption{Relative error in the $L_2$ norm for the approximation of the first $K=200$ coefficients of the POD decomposition (test set).}\label{bingml:errdom}
\end{center}
\end{table}

\subsection{Non-homogenious yield stress limit}
\begin{figure}[h!]
\includegraphics[width=0.9\textwidth]{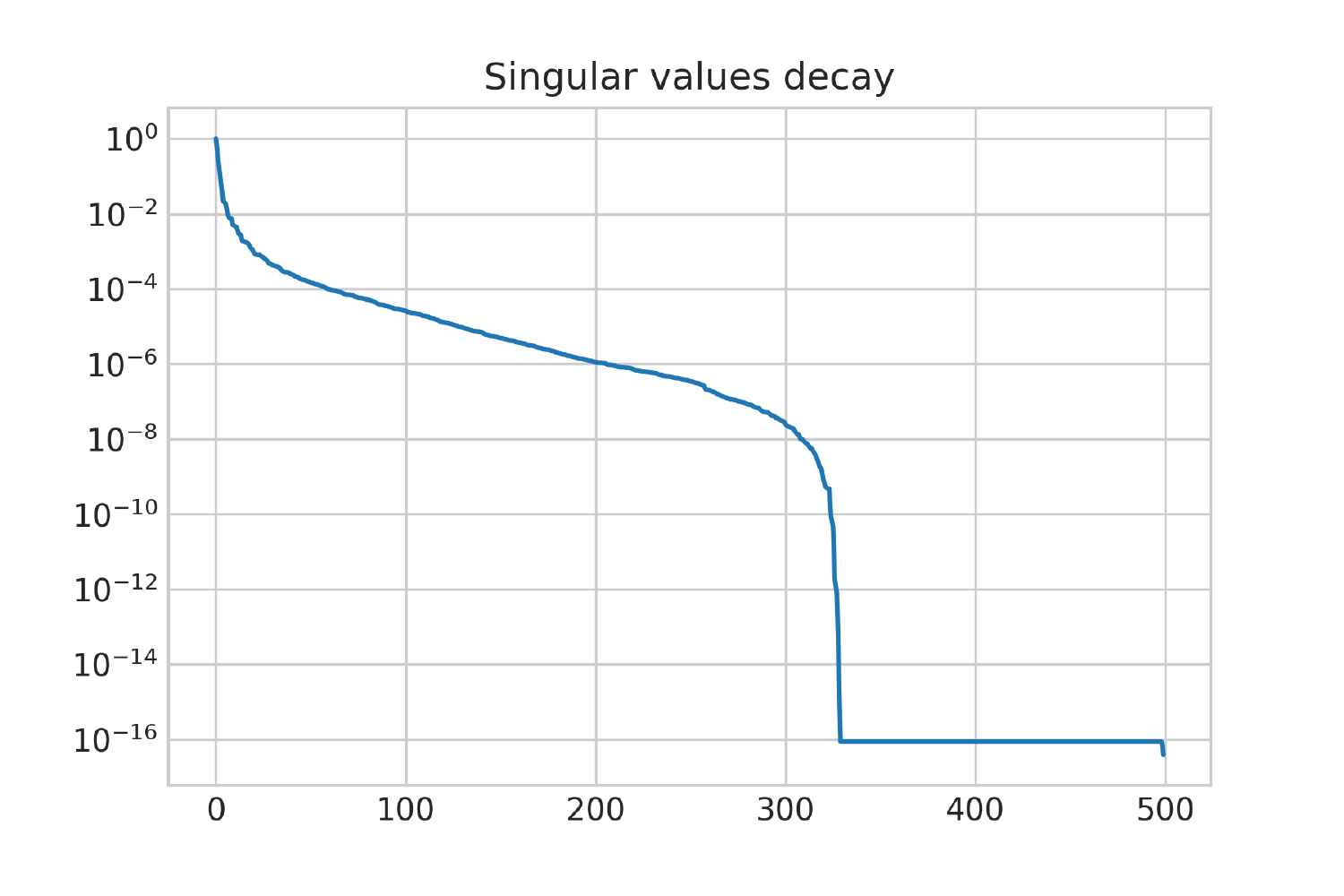}
\caption{Decay of singular values for two-parameter problem}\label{bingml:twoparam-decay}.
\end{figure}

\begin{figure}[h!]
\includegraphics[width=0.9\textwidth]{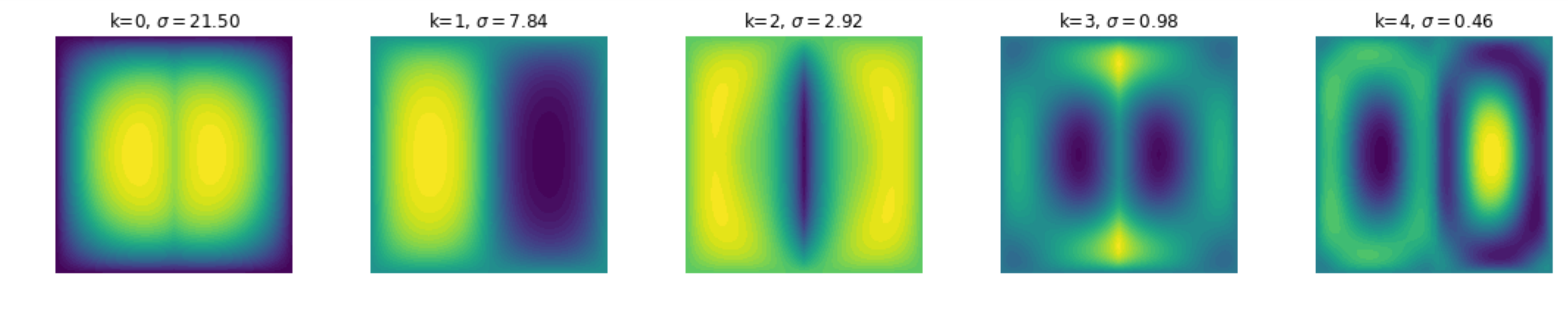}
\caption{First five POD basis functions}\label{bingml:twoparam-sv}.
\end{figure}

We consider the square case ($\Omega = [-\frac{1}{2}, \frac{1}{2}]^2$), and a piecewise-constant yield stress limit. This corresponds to the flow of different fluids with different yield stress limits. This type of flows is rather common in practice: it includes such processes as lamination, coextrusion and deposition. For the model problem
we consider a two-parameter problem with piecewise-constant $B(x_1, x_2)$:
$$
   B(x_1, x_2) = \begin{cases}
       B_1 & x_1 \leq 0, \\
       B_2 & x_1 > 0.
       \end{cases}
$$
The parameters $B_1$ and $B_2$ vary from $0$ to $0.8$. For this problem we need to select more snapshots than for a single stress limit to get similar accuracy. We also have to sample two-dimensional points in the parameter space.
To generate the POD basis, we compute $500$ snapshots with $B_1$ and $B_2$ generated using Halton quasi-random sequence generator \cite{niederreiter1992random}. The decay of singular values of the snapshot matrix is shown on Figure~\ref{bingml:twoparam-decay}, and several first singular vectors are depicted on Figure~\ref{bingml:twoparam-sv}.
Then we split the dataset randomly into training (67\% points) and testing (33\% points), and fit a fully connected ANN to map two input parameters ($B_1$ and $B_2$) to 
the first $20$  coefficients of the POD decomposition.
As an architecture of ANN we take $3$ hidden layers with $60, 50, 40$ hidden units, respectively and one output layer.

As an optimizer we again use the Adam optimizer with default parameters.
The relative error for the approximation on the test set is $0.008$. 

\begin{figure}[h!]
\includegraphics[width=0.9\textwidth]{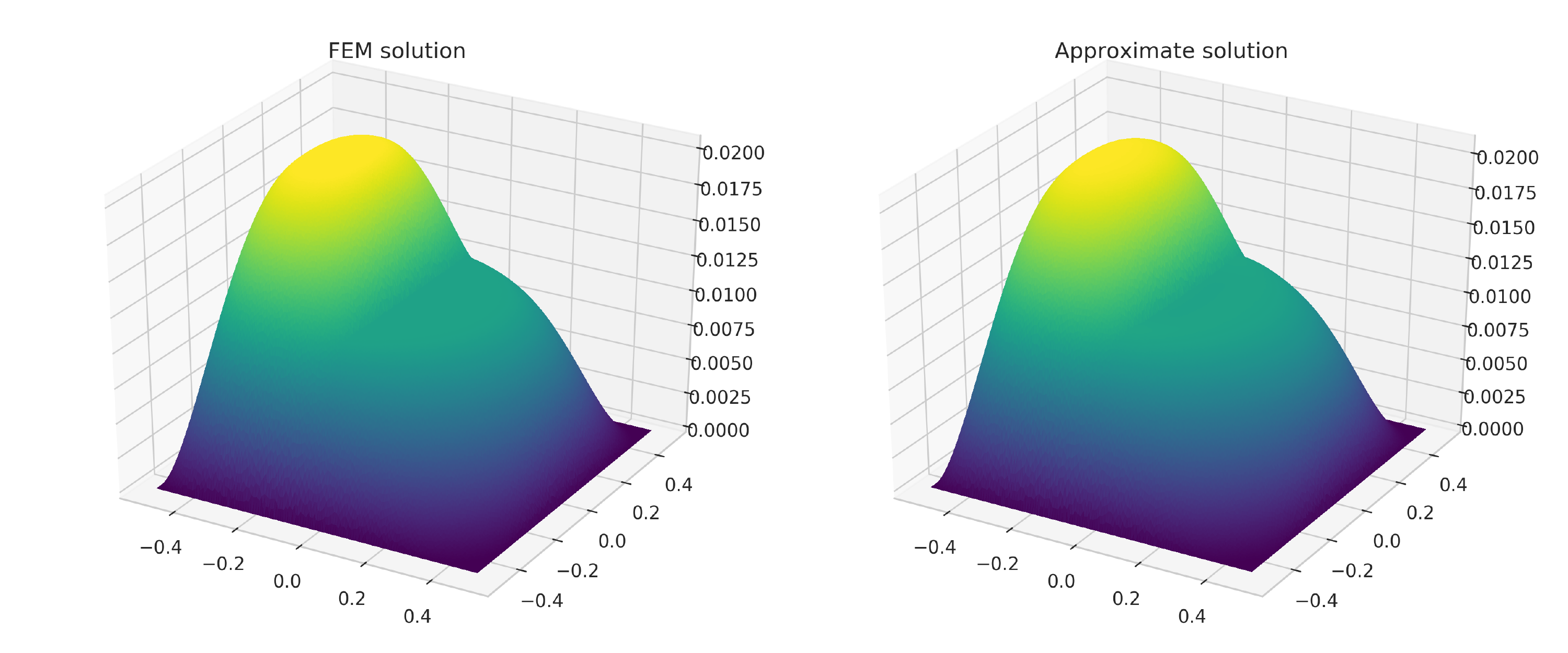}
\caption{Comparison of FEM and approximate solution}\label{bingml:twoparamcomp}.
\end{figure}

The main benefit of using the reduced model is that computing the approximate solution for new $B_1, B_2$ is very fast.
The final model takes approximately 4 milliseconds to evaluate, whereas the full FEM model takes 45 seconds. The comparison of two solutions is given on Figure~\ref{bingml:twoparamcomp} for $B_1 = 0.25$ and $B_2 = 0.05$.

\subsection{Flow rate vs pressure drop approximation}
\begin{figure}[h!]
  \includegraphics[width=0.9\textwidth]{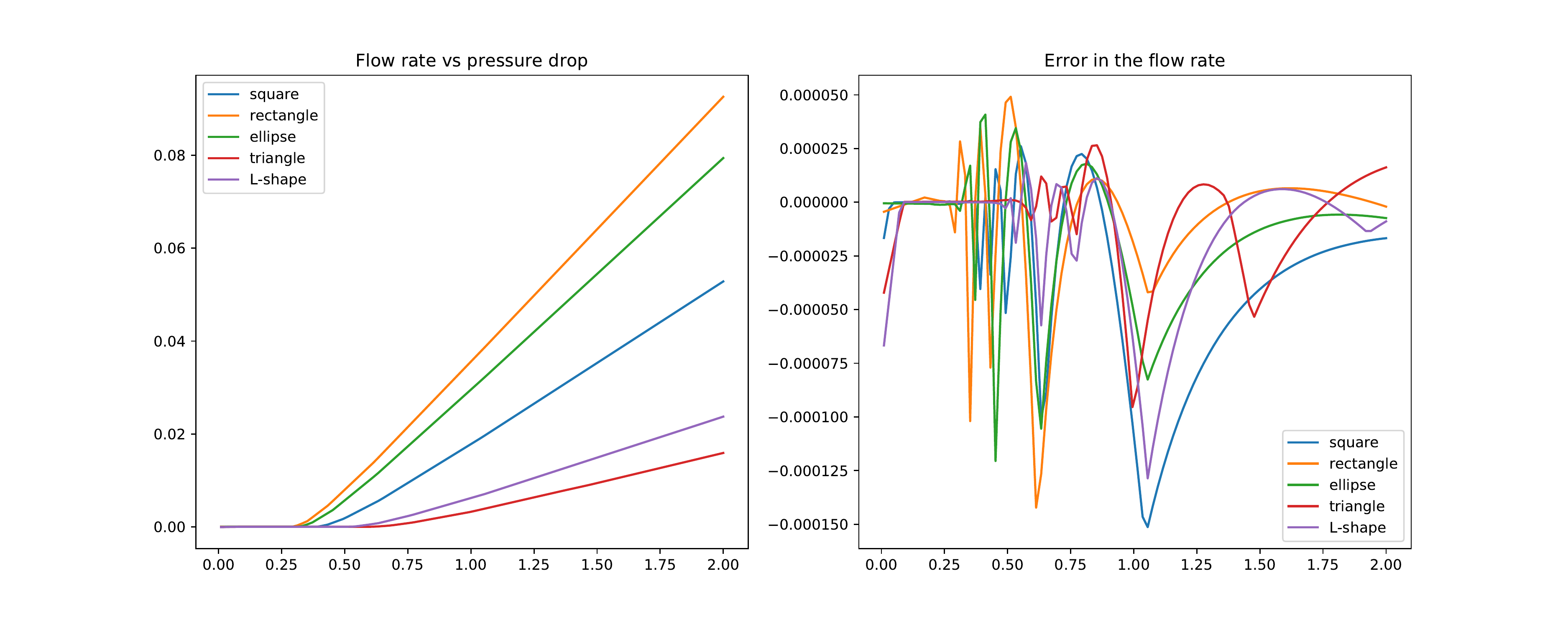}
  \caption{Flow rate vs pressure drop (left), difference between reduced and full model (right) for different cross-sections}\label{bingml:1par}
  \end{figure}

\begin{figure}[h!]
\centering
    \includegraphics[width=0.7\textwidth]{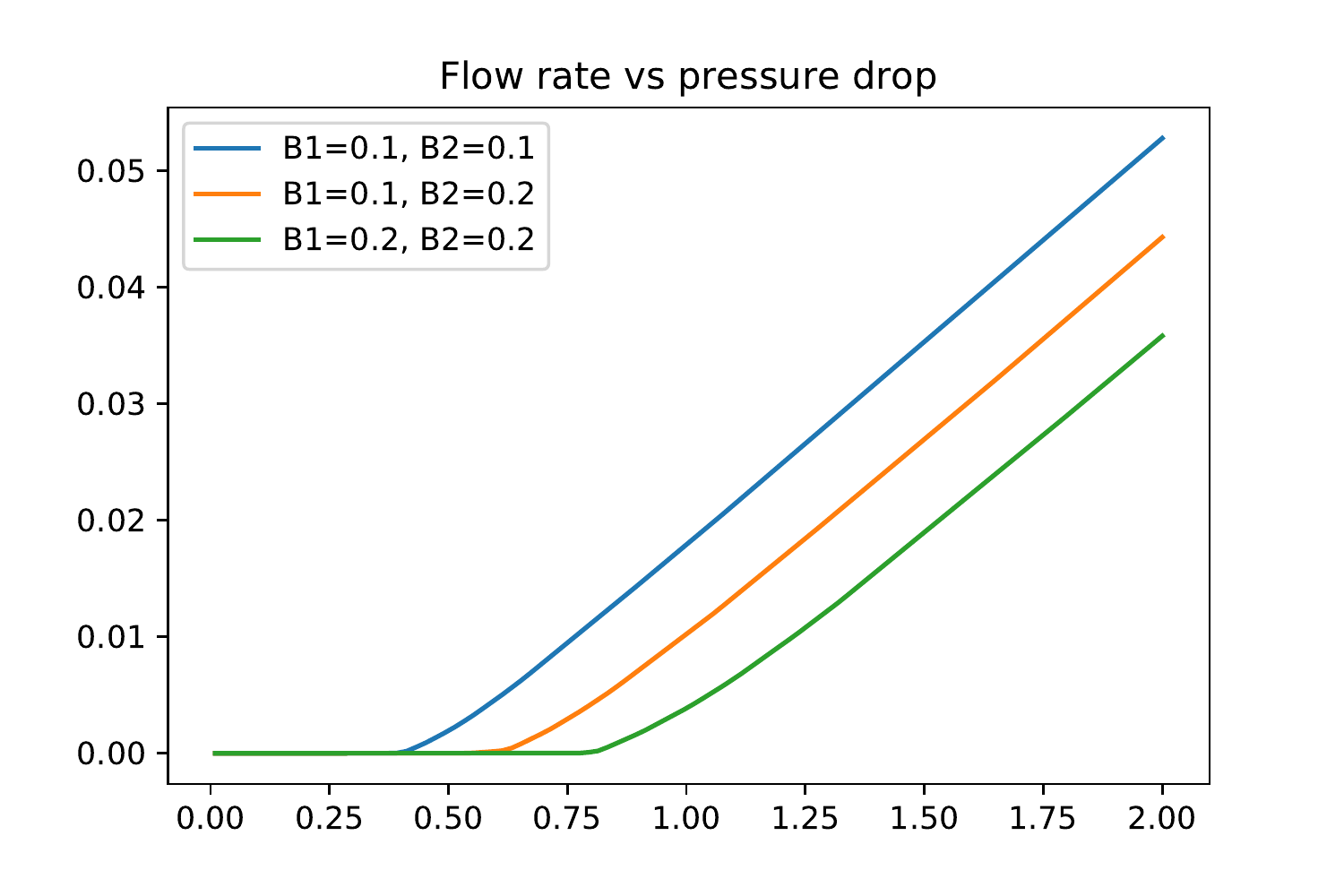}
\caption{Flow rate vs pressure drop for non-homogenious case for different Bingham numbers} \label{bingml:nh}.
\end{figure}

%To 
%If we have computed the approximant, we can use it to compute flow rate for any pressure drop very fast.
%Indeed, the solution with the pressure drop $\Delta p$ 
%and stress limit $B(x_1, x_2)$ can be computed from the solution with the pressure drop $11$ with stress limit $\frac{B(x_1, x_2)}{\Delta p}$ scaled by $\Delta p$.
As an example of application of the reduced order model, we use it to compute the dependence of the flow rate on the pressure drop. 
To do that, we need to solve \eqref{bingml:mosolov0} for fixed $\tau_s, \mu$ and varying pressure drop. This is equivalent 
to the computation of the solution of \eqref{bingml:mosolov} with varying $B$. 
We can also compare the result computed from the approximant and the solution obtained from the full model. 
On Figure ~\ref{bingml:1par}
the results are shown for different
cross-sections and $B=0.1$, and on Figure~\ref{bingml:nh} the results are shown for the non-homogenius yield stress case with different $B_1$ and $B_2$. The error behaviour is very similar
to one-parameter case, so we omit it here. 

It can be seen that in all cases the flow rate is approximated by the reduced model very accurately. The reduced order model is also much faster:
the total computation time for $100$ different pressure drops was $0.05$ seconds for the reduced model, and $66$ seconds for the full model.

\section{Conclusions and future work}
We have presented a general approach for the construction of reduced models of Bingham fluid flows
in the simplest possible case --- duct flow. Although being model, these flows share the main characteristics
of more general cases. The proposed method is ``easy to implement'': all the steps are automatic (generation of snapshots and fitting a neural network), but difficult to analyze: there is no guarantee that the approximated solution will share important properties
of the original solution, such as positivity. This requires a separate study, as well as the application of the constructed
reduced-order model to real-life design problems.

%\section{Acknowledgements}
%This work was supported by Russian Science Foundation grant 17-11-01376.
%\bibliographystyle{alpha}
%\bibliographystyle{apalike}
\bibliography{sample}

\begin{thebibliography}{35}%
\makeatletter
\providecommand \@ifxundefined [1]{%
 \@ifx{#1\undefined}
}%
\providecommand \@ifnum [1]{%
 \ifnum #1\expandafter \@firstoftwo
 \else \expandafter \@secondoftwo
 \fi
}%
\providecommand \@ifx [1]{%
 \ifx #1\expandafter \@firstoftwo
 \else \expandafter \@secondoftwo
 \fi
}%
\providecommand \natexlab [1]{#1}%
\providecommand \enquote  [1]{``#1''}%
\providecommand \bibnamefont  [1]{#1}%
\providecommand \bibfnamefont [1]{#1}%
\providecommand \citenamefont [1]{#1}%
\providecommand \href@noop [0]{\@secondoftwo}%
\providecommand \href [0]{\begingroup \@sanitize@url \@href}%
\providecommand \@href[1]{\@@startlink{#1}\@@href}%
\providecommand \@@href[1]{\endgroup#1\@@endlink}%
\providecommand \@sanitize@url [0]{\catcode `\\12\catcode `\$12\catcode
  `\&12\catcode `\#12\catcode `\^12\catcode `\_12\catcode `\%12\relax}%
\providecommand \@@startlink[1]{}%
\providecommand \@@endlink[0]{}%
\providecommand \url  [0]{\begingroup\@sanitize@url \@url }%
\providecommand \@url [1]{\endgroup\@href {#1}{\urlprefix }}%
\providecommand \urlprefix  [0]{URL }%
\providecommand \Eprint [0]{\href }%
\providecommand \doibase [0]{http://dx.doi.org/}%
\providecommand \selectlanguage [0]{\@gobble}%
\providecommand \bibinfo  [0]{\@secondoftwo}%
\providecommand \bibfield  [0]{\@secondoftwo}%
\providecommand \translation [1]{[#1]}%
\providecommand \BibitemOpen [0]{}%
\providecommand \bibitemStop [0]{}%
\providecommand \bibitemNoStop [0]{.\EOS\space}%
\providecommand \EOS [0]{\spacefactor3000\relax}%
\providecommand \BibitemShut  [1]{\csname bibitem#1\endcsname}%
\let\auto@bib@innerbib\@empty
%</preamble>
\bibitem [{\citenamefont {Frigaard}, \citenamefont {Paso},\ and\ \citenamefont
  {de~Souza~Mendes}(2017)}]{FrigPasoMend2017}%
  \BibitemOpen
  \bibfield  {author} {\bibinfo {author} {\bibfnamefont {I.~A.}\ \bibnamefont
  {Frigaard}}, \bibinfo {author} {\bibfnamefont {K.}~\bibnamefont {Paso}}, \
  and\ \bibinfo {author} {\bibfnamefont {P.}~\bibnamefont {de~Souza~Mendes}},\
  }\bibfield  {title} {\enquote {\bibinfo {title} {Bingham’s model in the oil
  and gas industry},}\ }\href@noop {} {\bibfield  {journal} {\bibinfo
  {journal} {Rheologica Acta}\ }\textbf {\bibinfo {volume} {56}},\ \bibinfo
  {pages} {259--282} (\bibinfo {year} {2017})}\BibitemShut {NoStop}%
\bibitem [{\citenamefont {Osiptsov}(2017)}]{osiptsov2017fluid}%
  \BibitemOpen
  \bibfield  {author} {\bibinfo {author} {\bibfnamefont {A.}~\bibnamefont
  {Osiptsov}},\ }\bibfield  {title} {\enquote {\bibinfo {title} {Fluid
  mechanics of hydraulic fracturing: A review},}\ }\href@noop {} {\bibfield
  {journal} {\bibinfo  {journal} {Journal of Petroleum Science and
  Engineering}\ } (\bibinfo {year} {2017})}\BibitemShut {NoStop}%
\bibitem [{\citenamefont {Guo}\ \emph {et~al.}(2017)\citenamefont {Guo},
  \citenamefont {Qu}, \citenamefont {Gong},\ and\ \citenamefont
  {Wang}}]{guo2017numerical}%
  \BibitemOpen
  \bibfield  {author} {\bibinfo {author} {\bibfnamefont {T.}~\bibnamefont
  {Guo}}, \bibinfo {author} {\bibfnamefont {Z.}~\bibnamefont {Qu}}, \bibinfo
  {author} {\bibfnamefont {F.}~\bibnamefont {Gong}}, \ and\ \bibinfo {author}
  {\bibfnamefont {X.}~\bibnamefont {Wang}},\ }\bibfield  {title} {\enquote
  {\bibinfo {title} {Numerical simulation of hydraulic fracture propagation
  guided by single radial boreholes},}\ }\href@noop {} {\bibfield  {journal}
  {\bibinfo  {journal} {Energies}\ }\textbf {\bibinfo {volume} {10}},\ \bibinfo
  {pages} {1680} (\bibinfo {year} {2017})}\BibitemShut {NoStop}%
\bibitem [{\citenamefont {Shojaei}, \citenamefont {Taleghani},\ and\
  \citenamefont {Li}(2014)}]{shojaei2014continuum}%
  \BibitemOpen
  \bibfield  {author} {\bibinfo {author} {\bibfnamefont {A.}~\bibnamefont
  {Shojaei}}, \bibinfo {author} {\bibfnamefont {A.}~\bibnamefont {Taleghani}},
  \ and\ \bibinfo {author} {\bibfnamefont {G.}~\bibnamefont {Li}},\ }\bibfield
  {title} {\enquote {\bibinfo {title} {A continuum damage failure model for
  hydraulic fracturing of porous rocks},}\ }\href@noop {} {\bibfield  {journal}
  {\bibinfo  {journal} {International Journal of Plasticity}\ }\textbf
  {\bibinfo {volume} {59}},\ \bibinfo {pages} {199--212} (\bibinfo {year}
  {2014})}\BibitemShut {NoStop}%
\bibitem [{\citenamefont {Bunger}\ and\ \citenamefont
  {Lecampion}(2017)}]{bunger2017four}%
  \BibitemOpen
  \bibfield  {author} {\bibinfo {author} {\bibfnamefont {A.}~\bibnamefont
  {Bunger}}\ and\ \bibinfo {author} {\bibfnamefont {B.}~\bibnamefont
  {Lecampion}},\ }\href@noop {} {\enquote {\bibinfo {title} {Four critical
  issues for successful hydraulic fracturing applications},}\ }\bibinfo {type}
  {Tech. Rep.}\ (\bibinfo  {institution} {CRC Press},\ \bibinfo {year}
  {2017})\BibitemShut {NoStop}%
\bibitem [{\citenamefont {Li}\ and\ \citenamefont
  {Shojaei}(2012)}]{li2012viscoplastic}%
  \BibitemOpen
  \bibfield  {author} {\bibinfo {author} {\bibfnamefont {G.}~\bibnamefont
  {Li}}\ and\ \bibinfo {author} {\bibfnamefont {A.}~\bibnamefont {Shojaei}},\
  }\bibfield  {title} {\enquote {\bibinfo {title} {A viscoplastic theory of
  shape memory polymer fibres with application to self-healing materials},}\
  }\href@noop {} {\bibfield  {journal} {\bibinfo  {journal} {Proc. R. Soc. A}\
  }\textbf {\bibinfo {volume} {468}},\ \bibinfo {pages} {2319--2346} (\bibinfo
  {year} {2012})}\BibitemShut {NoStop}%
\bibitem [{\citenamefont {Saramito}\ and\ \citenamefont
  {Wachs}(2017)}]{SarWachs2017}%
  \BibitemOpen
  \bibfield  {author} {\bibinfo {author} {\bibfnamefont {P.}~\bibnamefont
  {Saramito}}\ and\ \bibinfo {author} {\bibfnamefont {A.}~\bibnamefont
  {Wachs}},\ }\bibfield  {title} {\enquote {\bibinfo {title} {Progress in
  numerical simulation of yield stress fluid flows},}\ }\href@noop {}
  {\bibfield  {journal} {\bibinfo  {journal} {Rheologica Acta}\ }\textbf
  {\bibinfo {volume} {56}},\ \bibinfo {pages} {211--230} (\bibinfo {year}
  {2017})}\BibitemShut {NoStop}%
\bibitem [{\citenamefont {Mitsoulis}\ and\ \citenamefont
  {Tsamopoulos}(2017)}]{MitTsam2017}%
  \BibitemOpen
  \bibfield  {author} {\bibinfo {author} {\bibfnamefont {E.}~\bibnamefont
  {Mitsoulis}}\ and\ \bibinfo {author} {\bibfnamefont {J.}~\bibnamefont
  {Tsamopoulos}},\ }\bibfield  {title} {\enquote {\bibinfo {title} {Numerical
  simulations of complex yield-stress fluid flows},}\ }\href@noop {} {\bibfield
   {journal} {\bibinfo  {journal} {Rheologica Acta}\ }\textbf {\bibinfo
  {volume} {56}},\ \bibinfo {pages} {231--258} (\bibinfo {year}
  {2017})}\BibitemShut {NoStop}%
\bibitem [{\citenamefont {Ammar}\ \emph {et~al.}(2006)\citenamefont {Ammar},
  \citenamefont {Mokdad}, \citenamefont {Chinesta},\ and\ \citenamefont
  {Keunings}}]{ammar2006new}%
  \BibitemOpen
  \bibfield  {author} {\bibinfo {author} {\bibfnamefont {A.}~\bibnamefont
  {Ammar}}, \bibinfo {author} {\bibfnamefont {B.}~\bibnamefont {Mokdad}},
  \bibinfo {author} {\bibfnamefont {F.}~\bibnamefont {Chinesta}}, \ and\
  \bibinfo {author} {\bibfnamefont {R.}~\bibnamefont {Keunings}},\ }\bibfield
  {title} {\enquote {\bibinfo {title} {A new family of solvers for some classes
  of multidimensional partial differential equations encountered in kinetic
  theory modeling of complex fluids},}\ }\href@noop {} {\bibfield  {journal}
  {\bibinfo  {journal} {Journal of Non-Newtonian Fluid Mechanics}\ }\textbf
  {\bibinfo {volume} {139}},\ \bibinfo {pages} {153--176} (\bibinfo {year}
  {2006})}\BibitemShut {NoStop}%
\bibitem [{\citenamefont {Chinesta}\ \emph {et~al.}(2011)\citenamefont
  {Chinesta}, \citenamefont {Ammar}, \citenamefont {Leygue},\ and\
  \citenamefont {Keunings}}]{chinesta2011overview}%
  \BibitemOpen
  \bibfield  {author} {\bibinfo {author} {\bibfnamefont {F.}~\bibnamefont
  {Chinesta}}, \bibinfo {author} {\bibfnamefont {A.}~\bibnamefont {Ammar}},
  \bibinfo {author} {\bibfnamefont {A.}~\bibnamefont {Leygue}}, \ and\ \bibinfo
  {author} {\bibfnamefont {R.}~\bibnamefont {Keunings}},\ }\bibfield  {title}
  {\enquote {\bibinfo {title} {An overview of the proper generalized
  decomposition with applications in computational rheology},}\ }\href@noop {}
  {\bibfield  {journal} {\bibinfo  {journal} {Journal of Non-Newtonian Fluid
  Mechanics}\ }\textbf {\bibinfo {volume} {166}},\ \bibinfo {pages} {578--592}
  (\bibinfo {year} {2011})}\BibitemShut {NoStop}%
\bibitem [{\citenamefont {Mosolov}\ and\ \citenamefont
  {Miasnikov}(1965)}]{MosMyas1965}%
  \BibitemOpen
  \bibfield  {author} {\bibinfo {author} {\bibfnamefont {P.}~\bibnamefont
  {Mosolov}}\ and\ \bibinfo {author} {\bibfnamefont {V.}~\bibnamefont
  {Miasnikov}},\ }\bibfield  {title} {\enquote {\bibinfo {title} {Variational
  methods in the theory of the fluidity of a viscous-plastic medium},}\
  }\href@noop {} {\bibfield  {journal} {\bibinfo  {journal} {Journal of Applied
  Mathematics and Mechanics}\ }\textbf {\bibinfo {volume} {29}},\ \bibinfo
  {pages} {545--577} (\bibinfo {year} {1965})}\BibitemShut {NoStop}%
\bibitem [{\citenamefont {Mosolov}\ and\ \citenamefont
  {Miasnikov}(1966)}]{MosMyas1966}%
  \BibitemOpen
  \bibfield  {author} {\bibinfo {author} {\bibfnamefont {P.}~\bibnamefont
  {Mosolov}}\ and\ \bibinfo {author} {\bibfnamefont {V.}~\bibnamefont
  {Miasnikov}},\ }\bibfield  {title} {\enquote {\bibinfo {title} {On stagnant
  flow regions of a viscous-plastic medium in pipes},}\ }\href@noop {}
  {\bibfield  {journal} {\bibinfo  {journal} {Journal of Applied Mathematics
  and Mechanics}\ }\textbf {\bibinfo {volume} {30}},\ \bibinfo {pages}
  {841--854} (\bibinfo {year} {1966})}\BibitemShut {NoStop}%
\bibitem [{\citenamefont {Mosolov}\ and\ \citenamefont
  {Miasnikov}(1967)}]{MosMyas1967}%
  \BibitemOpen
  \bibfield  {author} {\bibinfo {author} {\bibfnamefont {P.}~\bibnamefont
  {Mosolov}}\ and\ \bibinfo {author} {\bibfnamefont {V.}~\bibnamefont
  {Miasnikov}},\ }\bibfield  {title} {\enquote {\bibinfo {title} {On
  qualitative singularities of the flow of a viscoplastic medium in pipes},}\
  }\href@noop {} {\bibfield  {journal} {\bibinfo  {journal} {Journal of Applied
  Mathematics and Mechanics}\ }\textbf {\bibinfo {volume} {31}},\ \bibinfo
  {pages} {609--613} (\bibinfo {year} {1967})}\BibitemShut {NoStop}%
\bibitem [{\citenamefont {Huilgol}\ and\ \citenamefont
  {You}(2005)}]{Huilgol2005}%
  \BibitemOpen
  \bibfield  {author} {\bibinfo {author} {\bibfnamefont {R.}~\bibnamefont
  {Huilgol}}\ and\ \bibinfo {author} {\bibfnamefont {Z.}~\bibnamefont {You}},\
  }\bibfield  {title} {\enquote {\bibinfo {title} {Application of the augmented
  {Lagrangian} method to steady pipe flows of {Bingham}, {Casson} and
  {Herschel–Bulkley} fluids},}\ }\href@noop {} {\bibfield  {journal}
  {\bibinfo  {journal} {Journal of Non-Newtonian Fluid Mechanics}\ }\textbf
  {\bibinfo {volume} {128}},\ \bibinfo {pages} {126--143} (\bibinfo {year}
  {2005})}\BibitemShut {NoStop}%
\bibitem [{\citenamefont {Saramito}\ and\ \citenamefont
  {Roquet}(2001)}]{SaramRoquet2001}%
  \BibitemOpen
  \bibfield  {author} {\bibinfo {author} {\bibfnamefont {P.}~\bibnamefont
  {Saramito}}\ and\ \bibinfo {author} {\bibfnamefont {N.}~\bibnamefont
  {Roquet}},\ }\bibfield  {title} {\enquote {\bibinfo {title} {An adaptive
  finite element method for viscoplastic fluid flows in pipes},}\ }\href@noop
  {} {\bibfield  {journal} {\bibinfo  {journal} {Computer Methods in Applied
  Mechanics and Engineering}\ }\textbf {\bibinfo {volume} {190}},\ \bibinfo
  {pages} {5391--5412} (\bibinfo {year} {2001})}\BibitemShut {NoStop}%
\bibitem [{\citenamefont {Moyers-Gonzalez}\ and\ \citenamefont
  {A.}(2004)}]{MoyFrig2004}%
  \BibitemOpen
  \bibfield  {author} {\bibinfo {author} {\bibfnamefont {M.~A.}\ \bibnamefont
  {Moyers-Gonzalez}}\ and\ \bibinfo {author} {\bibfnamefont {F.~I.}\
  \bibnamefont {A.}},\ }\bibfield  {title} {\enquote {\bibinfo {title}
  {Numerical solution of duct flows of multiple visco-plastic fluids},}\
  }\href@noop {} {\bibfield  {journal} {\bibinfo  {journal} {Journal of
  Non-Newtonian {Fluid} Mechanics}\ }\textbf {\bibinfo {volume} {122}},\
  \bibinfo {pages} {227--241} (\bibinfo {year} {2004})}\BibitemShut {NoStop}%
\bibitem [{\citenamefont {Muravleva}(2009)}]{MuravMM2009}%
  \BibitemOpen
  \bibfield  {author} {\bibinfo {author} {\bibfnamefont {E.}~\bibnamefont
  {Muravleva}},\ }\bibfield  {title} {\enquote {\bibinfo {title}
  {Finite-difference schemes for the computation of viscoplastic medium flows
  in a channel},}\ }\href@noop {} {\bibfield  {journal} {\bibinfo  {journal}
  {Mathematical Models and Computer Simulations}\ }\textbf {\bibinfo {volume}
  {1}},\ \bibinfo {pages} {768--279} (\bibinfo {year} {2009})}\BibitemShut
  {NoStop}%
\bibitem [{\citenamefont {Walton}\ and\ \citenamefont
  {Bittleston}(1991)}]{walton1991axial}%
  \BibitemOpen
  \bibfield  {author} {\bibinfo {author} {\bibfnamefont {I.~C.}\ \bibnamefont
  {Walton}}\ and\ \bibinfo {author} {\bibfnamefont {S.~H.}\ \bibnamefont
  {Bittleston}},\ }\bibfield  {title} {\enquote {\bibinfo {title} {The axial
  flow of a {Bingham} plastic in a narrow eccentric annulus},}\ }\href@noop {}
  {\bibfield  {journal} {\bibinfo  {journal} {Journal of Fluid Mechanics}\
  }\textbf {\bibinfo {volume} {222}},\ \bibinfo {pages} {39--60} (\bibinfo
  {year} {1991})}\BibitemShut {NoStop}%
\bibitem [{\citenamefont {Szabo}\ and\ \citenamefont
  {Hassager}(1992)}]{szabo1992flow}%
  \BibitemOpen
  \bibfield  {author} {\bibinfo {author} {\bibfnamefont {P.}~\bibnamefont
  {Szabo}}\ and\ \bibinfo {author} {\bibfnamefont {O.}~\bibnamefont
  {Hassager}},\ }\bibfield  {title} {\enquote {\bibinfo {title} {Flow of
  viscoplastic fluids in eccentric annular geometries},}\ }\href@noop {}
  {\bibfield  {journal} {\bibinfo  {journal} {Journal of Non-Newtonian Fluid
  Mechanics}\ }\textbf {\bibinfo {volume} {45}},\ \bibinfo {pages} {149--169}
  (\bibinfo {year} {1992})}\BibitemShut {NoStop}%
\bibitem [{\citenamefont {Wachs}(2007)}]{wachs2007numerical}%
  \BibitemOpen
  \bibfield  {author} {\bibinfo {author} {\bibfnamefont {A.}~\bibnamefont
  {Wachs}},\ }\bibfield  {title} {\enquote {\bibinfo {title} {Numerical
  simulation of steady {Bingham} flow through an eccentric annular
  cross-section by distributed {Lagrange} multiplier/fictitious domain and
  augmented {Lagrangian} methods},}\ }\href@noop {} {\bibfield  {journal}
  {\bibinfo  {journal} {Journal of {Non-Newtonian} Fluid Mechanics}\ }\textbf
  {\bibinfo {volume} {142}},\ \bibinfo {pages} {183--198} (\bibinfo {year}
  {2007})}\BibitemShut {NoStop}%
\bibitem [{\citenamefont {Muravleva}\ and\ \citenamefont
  {Muravleva}(2009)}]{muravleva2009unsteady}%
  \BibitemOpen
  \bibfield  {author} {\bibinfo {author} {\bibfnamefont {E.~A.}\ \bibnamefont
  {Muravleva}}\ and\ \bibinfo {author} {\bibfnamefont {L.~V.}\ \bibnamefont
  {Muravleva}},\ }\bibfield  {title} {\enquote {\bibinfo {title} {Unsteady
  flows of a viscoplastic medium in channels},}\ }\href@noop {} {\bibfield
  {journal} {\bibinfo  {journal} {Mechanics of Solids}\ }\textbf {\bibinfo
  {volume} {44}},\ \bibinfo {pages} {792--812} (\bibinfo {year}
  {2009})}\BibitemShut {NoStop}%
\bibitem [{\citenamefont {Muravleva}\ \emph {et~al.}(2010)\citenamefont
  {Muravleva}, \citenamefont {Muravleva}, \citenamefont {Georgiou},\ and\
  \citenamefont {Mitsoulis}}]{Mur_JNNFM2010}%
  \BibitemOpen
  \bibfield  {author} {\bibinfo {author} {\bibfnamefont {L.}~\bibnamefont
  {Muravleva}}, \bibinfo {author} {\bibfnamefont {E.}~\bibnamefont
  {Muravleva}}, \bibinfo {author} {\bibfnamefont {G.}~\bibnamefont {Georgiou}},
  \ and\ \bibinfo {author} {\bibfnamefont {E.}~\bibnamefont {Mitsoulis}},\
  }\bibfield  {title} {\enquote {\bibinfo {title} {Numerical simulations of
  cessation flows of a {Bingham} plastic with the augmented {Lagrangian}
  method},}\ }\href@noop {} {\bibfield  {journal} {\bibinfo  {journal} {Journal
  of Non-Newtonian Fluid Mechanics}\ }\textbf {\bibinfo {volume} {164}},\
  \bibinfo {pages} {544--550} (\bibinfo {year} {2010})}\BibitemShut {NoStop}%
\bibitem [{\citenamefont {Roquet}\ and\ \citenamefont
  {Saramito}(2008)}]{roquet2008adaptive}%
  \BibitemOpen
  \bibfield  {author} {\bibinfo {author} {\bibfnamefont {N.}~\bibnamefont
  {Roquet}}\ and\ \bibinfo {author} {\bibfnamefont {P.}~\bibnamefont
  {Saramito}},\ }\bibfield  {title} {\enquote {\bibinfo {title} {An adaptive
  finite element method for viscoplastic flows in a square pipe with
  stick--slip at the wall},}\ }\href@noop {} {\bibfield  {journal} {\bibinfo
  {journal} {Journal of Non-Newtonian Fluid Mechanics}\ }\textbf {\bibinfo
  {volume} {155}},\ \bibinfo {pages} {101--115} (\bibinfo {year}
  {2008})}\BibitemShut {NoStop}%
\bibitem [{\citenamefont {Damianou}\ \emph {et~al.}(2014)\citenamefont
  {Damianou}, \citenamefont {Philippou}, \citenamefont {Kaoullas},\ and\
  \citenamefont {Georgiou}}]{DamPhil}%
  \BibitemOpen
  \bibfield  {author} {\bibinfo {author} {\bibfnamefont {Y.}~\bibnamefont
  {Damianou}}, \bibinfo {author} {\bibfnamefont {M.}~\bibnamefont {Philippou}},
  \bibinfo {author} {\bibfnamefont {G.}~\bibnamefont {Kaoullas}}, \ and\
  \bibinfo {author} {\bibfnamefont {G.}~\bibnamefont {Georgiou}},\ }\bibfield
  {title} {\enquote {\bibinfo {title} {Cessation of viscoplastic {Poiseuille}
  flow with wall slip},}\ }\href@noop {} {\bibfield  {journal} {\bibinfo
  {journal} {Journal of Non-Newtonian Fluid Mechanics}\ }\textbf {\bibinfo
  {volume} {203}},\ \bibinfo {pages} {24--37} (\bibinfo {year}
  {2014})}\BibitemShut {NoStop}%
\bibitem [{\citenamefont {Damianou}\ and\ \citenamefont
  {Georgiou}(2014)}]{DamGeor}%
  \BibitemOpen
  \bibfield  {author} {\bibinfo {author} {\bibfnamefont {Y.}~\bibnamefont
  {Damianou}}\ and\ \bibinfo {author} {\bibfnamefont {G.}~\bibnamefont
  {Georgiou}},\ }\bibfield  {title} {\enquote {\bibinfo {title} {Viscoplastic
  {Poiseuille} flow in a rectangular duct with wall slip},}\ }\href@noop {}
  {\bibfield  {journal} {\bibinfo  {journal} {Journal of Non-Newtonian Fluid
  Mechanics}\ }\textbf {\bibinfo {volume} {214}},\ \bibinfo {pages} {88--105}
  (\bibinfo {year} {2014})}\BibitemShut {NoStop}%
\bibitem [{\citenamefont {Fortin}\ and\ \citenamefont
  {Glowinski}(1983)}]{Glowinski}%
  \BibitemOpen
  \bibfield  {author} {\bibinfo {author} {\bibfnamefont {M.}~\bibnamefont
  {Fortin}}\ and\ \bibinfo {author} {\bibfnamefont {R.}~\bibnamefont
  {Glowinski}},\ }\href@noop {} {\enquote {\bibinfo {title} {The augmented
  {Lagrangian} method},}\ } (\bibinfo {year} {1983})\BibitemShut {NoStop}%
\bibitem [{\citenamefont {Logg}, \citenamefont {Mardal},\ and\ \citenamefont
  {Wells}(2012)}]{LoggMardalEtAl2012a}%
  \BibitemOpen
  \bibfield  {author} {\bibinfo {author} {\bibfnamefont {A.}~\bibnamefont
  {Logg}}, \bibinfo {author} {\bibfnamefont {K.}~\bibnamefont {Mardal}}, \ and\
  \bibinfo {author} {\bibfnamefont {G.}~\bibnamefont {Wells}},\ }\href
  {\doibase 10.1007/978-3-642-23099-8} {\emph {\bibinfo {title} {Automated
  Solution of Differential Equations by the Finite Element Method}}}\ (\bibinfo
   {publisher} {Springer},\ \bibinfo {year} {2012})\BibitemShut {NoStop}%
\bibitem [{\citenamefont {Kunisch}\ and\ \citenamefont
  {Volkwein}(2002)}]{kunisch2002galerkin}%
  \BibitemOpen
  \bibfield  {author} {\bibinfo {author} {\bibfnamefont {K.}~\bibnamefont
  {Kunisch}}\ and\ \bibinfo {author} {\bibfnamefont {S.}~\bibnamefont
  {Volkwein}},\ }\bibfield  {title} {\enquote {\bibinfo {title} {Galerkin
  proper orthogonal decomposition methods for a general equation in fluid
  dynamics},}\ }\href@noop {} {\bibfield  {journal} {\bibinfo  {journal} {SIAM
  Journal on Numerical analysis}\ }\textbf {\bibinfo {volume} {40}},\ \bibinfo
  {pages} {492--515} (\bibinfo {year} {2002})}\BibitemShut {NoStop}%
\bibitem [{\citenamefont {Peterson}(1989)}]{peterson1989reduced}%
  \BibitemOpen
  \bibfield  {author} {\bibinfo {author} {\bibfnamefont {J.~S.}\ \bibnamefont
  {Peterson}},\ }\bibfield  {title} {\enquote {\bibinfo {title} {The reduced
  basis method for incompressible viscous flow calculations},}\ }\href@noop {}
  {\bibfield  {journal} {\bibinfo  {journal} {SIAM Journal on Scientific and
  Statistical Computing}\ }\textbf {\bibinfo {volume} {10}},\ \bibinfo {pages}
  {777--786} (\bibinfo {year} {1989})}\BibitemShut {NoStop}%
\bibitem [{\citenamefont {Lassila}\ \emph {et~al.}(2014)\citenamefont
  {Lassila}, \citenamefont {Manzoni}, \citenamefont {Quarteroni},\ and\
  \citenamefont {Rozza}}]{lassila2014model}%
  \BibitemOpen
  \bibfield  {author} {\bibinfo {author} {\bibfnamefont {T.}~\bibnamefont
  {Lassila}}, \bibinfo {author} {\bibfnamefont {A.}~\bibnamefont {Manzoni}},
  \bibinfo {author} {\bibfnamefont {A.}~\bibnamefont {Quarteroni}}, \ and\
  \bibinfo {author} {\bibfnamefont {G.}~\bibnamefont {Rozza}},\ }\bibfield
  {title} {\enquote {\bibinfo {title} {Model order reduction in fluid dynamics:
  challenges and perspectives},}\ }in\ \href@noop {} {\emph {\bibinfo
  {booktitle} {Reduced Order Methods for modeling and computational
  reduction}}}\ (\bibinfo  {publisher} {Springer},\ \bibinfo {year} {2014})\
  pp.\ \bibinfo {pages} {235--273}\BibitemShut {NoStop}%
\bibitem [{\citenamefont {Goodfellow}, \citenamefont {Bengio},\ and\
  \citenamefont {Courville}(2016)}]{Goodfellow-et-al-2016}%
  \BibitemOpen
  \bibfield  {author} {\bibinfo {author} {\bibfnamefont {I.}~\bibnamefont
  {Goodfellow}}, \bibinfo {author} {\bibfnamefont {Y.}~\bibnamefont {Bengio}},
  \ and\ \bibinfo {author} {\bibfnamefont {A.}~\bibnamefont {Courville}},\
  }\href@noop {} {\emph {\bibinfo {title} {Deep Learning}}}\ (\bibinfo
  {publisher} {MIT Press},\ \bibinfo {year} {2016})\ \bibinfo {note}
  {\url{http://www.deeplearningbook.org}}\BibitemShut {NoStop}%
\bibitem [{\citenamefont {Hornik}, \citenamefont {Stinchcombe},\ and\
  \citenamefont {White}(1989)}]{hornik1989multilayer}%
  \BibitemOpen
  \bibfield  {author} {\bibinfo {author} {\bibfnamefont {K.}~\bibnamefont
  {Hornik}}, \bibinfo {author} {\bibfnamefont {M.}~\bibnamefont {Stinchcombe}},
  \ and\ \bibinfo {author} {\bibfnamefont {H.}~\bibnamefont {White}},\
  }\bibfield  {title} {\enquote {\bibinfo {title} {Multilayer feedforward
  networks are universal approximators},}\ }\href@noop {} {\bibfield  {journal}
  {\bibinfo  {journal} {Neural networks}\ }\textbf {\bibinfo {volume} {2}},\
  \bibinfo {pages} {359--366} (\bibinfo {year} {1989})}\BibitemShut {NoStop}%
\bibitem [{\citenamefont {Abadi}\ \emph {et~al.}(2015)\citenamefont {Abadi},
  \citenamefont {Agarwal}, \citenamefont {Barham}, \citenamefont {Brevdo},
  \citenamefont {Chen}, \citenamefont {Citro}, \citenamefont {Corrado},
  \citenamefont {Davis}, \citenamefont {Dean}, \citenamefont {Devin},
  \citenamefont {Ghemawat}, \citenamefont {Goodfellow}, \citenamefont {Harp},
  \citenamefont {Irving}, \citenamefont {Isard}, \citenamefont {Jia},
  \citenamefont {Jozefowicz}, \citenamefont {Kaiser}, \citenamefont {Kudlur},
  \citenamefont {Levenberg}, \citenamefont {Man\'{e}}, \citenamefont {Monga},
  \citenamefont {Moore}, \citenamefont {Murray}, \citenamefont {Olah},
  \citenamefont {Schuster}, \citenamefont {Shlens}, \citenamefont {Steiner},
  \citenamefont {Sutskever}, \citenamefont {Talwar}, \citenamefont {Tucker},
  \citenamefont {Vanhoucke}, \citenamefont {Vasudevan}, \citenamefont
  {Vi\'{e}gas}, \citenamefont {Vinyals}, \citenamefont {Warden}, \citenamefont
  {Wattenberg}, \citenamefont {Wicke}, \citenamefont {Yu},\ and\ \citenamefont
  {Zheng}}]{tensorflow2015-whitepaper}%
  \BibitemOpen
  \bibfield  {author} {\bibinfo {author} {\bibfnamefont {M.}~\bibnamefont
  {Abadi}}, \bibinfo {author} {\bibfnamefont {A.}~\bibnamefont {Agarwal}},
  \bibinfo {author} {\bibfnamefont {P.}~\bibnamefont {Barham}}, \bibinfo
  {author} {\bibfnamefont {E.}~\bibnamefont {Brevdo}}, \bibinfo {author}
  {\bibfnamefont {Z.}~\bibnamefont {Chen}}, \bibinfo {author} {\bibfnamefont
  {C.}~\bibnamefont {Citro}}, \bibinfo {author} {\bibfnamefont {G.~S.}\
  \bibnamefont {Corrado}}, \bibinfo {author} {\bibfnamefont {A.}~\bibnamefont
  {Davis}}, \bibinfo {author} {\bibfnamefont {J.}~\bibnamefont {Dean}},
  \bibinfo {author} {\bibfnamefont {M.}~\bibnamefont {Devin}}, \bibinfo
  {author} {\bibfnamefont {S.}~\bibnamefont {Ghemawat}}, \bibinfo {author}
  {\bibfnamefont {I.}~\bibnamefont {Goodfellow}}, \bibinfo {author}
  {\bibfnamefont {A.}~\bibnamefont {Harp}}, \bibinfo {author} {\bibfnamefont
  {G.}~\bibnamefont {Irving}}, \bibinfo {author} {\bibfnamefont
  {M.}~\bibnamefont {Isard}}, \bibinfo {author} {\bibfnamefont
  {Y.}~\bibnamefont {Jia}}, \bibinfo {author} {\bibfnamefont {R.}~\bibnamefont
  {Jozefowicz}}, \bibinfo {author} {\bibfnamefont {L.}~\bibnamefont {Kaiser}},
  \bibinfo {author} {\bibfnamefont {M.}~\bibnamefont {Kudlur}}, \bibinfo
  {author} {\bibfnamefont {J.}~\bibnamefont {Levenberg}}, \bibinfo {author}
  {\bibfnamefont {D.}~\bibnamefont {Man\'{e}}}, \bibinfo {author}
  {\bibfnamefont {R.}~\bibnamefont {Monga}}, \bibinfo {author} {\bibfnamefont
  {S.}~\bibnamefont {Moore}}, \bibinfo {author} {\bibfnamefont
  {D.}~\bibnamefont {Murray}}, \bibinfo {author} {\bibfnamefont
  {C.}~\bibnamefont {Olah}}, \bibinfo {author} {\bibfnamefont {M.}~\bibnamefont
  {Schuster}}, \bibinfo {author} {\bibfnamefont {J.}~\bibnamefont {Shlens}},
  \bibinfo {author} {\bibfnamefont {B.}~\bibnamefont {Steiner}}, \bibinfo
  {author} {\bibfnamefont {I.}~\bibnamefont {Sutskever}}, \bibinfo {author}
  {\bibfnamefont {K.}~\bibnamefont {Talwar}}, \bibinfo {author} {\bibfnamefont
  {P.}~\bibnamefont {Tucker}}, \bibinfo {author} {\bibfnamefont
  {V.}~\bibnamefont {Vanhoucke}}, \bibinfo {author} {\bibfnamefont
  {V.}~\bibnamefont {Vasudevan}}, \bibinfo {author} {\bibfnamefont
  {F.}~\bibnamefont {Vi\'{e}gas}}, \bibinfo {author} {\bibfnamefont
  {O.}~\bibnamefont {Vinyals}}, \bibinfo {author} {\bibfnamefont
  {P.}~\bibnamefont {Warden}}, \bibinfo {author} {\bibfnamefont
  {M.}~\bibnamefont {Wattenberg}}, \bibinfo {author} {\bibfnamefont
  {M.}~\bibnamefont {Wicke}}, \bibinfo {author} {\bibfnamefont
  {Y.}~\bibnamefont {Yu}}, \ and\ \bibinfo {author} {\bibfnamefont
  {X.}~\bibnamefont {Zheng}},\ }\href {https://www.tensorflow.org/} {\enquote
  {\bibinfo {title} {{TensorFlow}: Large-scale machine learning on
  heterogeneous systems},}\ } (\bibinfo {year} {2015}),\ \bibinfo {note}
  {software available from tensorflow.org}\BibitemShut {NoStop}%
\bibitem [{\citenamefont {Kingma}\ and\ \citenamefont
  {Ba}(2014)}]{kingma2014adam}%
  \BibitemOpen
  \bibfield  {author} {\bibinfo {author} {\bibfnamefont {D.}~\bibnamefont
  {Kingma}}\ and\ \bibinfo {author} {\bibfnamefont {J.}~\bibnamefont {Ba}},\
  }\bibfield  {title} {\enquote {\bibinfo {title} {{ADAM}: A method for
  stochastic optimization},}\ }\href@noop {} {\bibfield  {journal} {\bibinfo
  {journal} {arXiv preprint arXiv:1412.6980}\ } (\bibinfo {year}
  {2014})}\BibitemShut {NoStop}%
\bibitem [{\citenamefont {Niederreiter}(1992)}]{niederreiter1992random}%
  \BibitemOpen
  \bibfield  {author} {\bibinfo {author} {\bibfnamefont {H.}~\bibnamefont
  {Niederreiter}},\ }\href@noop {} {\emph {\bibinfo {title} {Random number
  generation and quasi-{Monte} {Carlo} methods}}}\ (\bibinfo  {publisher}
  {SIAM},\ \bibinfo {year} {1992})\BibitemShut {NoStop}%
\end{thebibliography}%

\end{document}